\documentclass[11pt]{article}
\usepackage{amsmath}
\usepackage{amssymb}
\usepackage{graphics}
\usepackage{graphicx}
\usepackage{array}
\newtheorem{theo}{\hspace*{\parindent}Theorem}

\newtheorem{lemma}{\hspace*{\parindent}Lemma}

\def\arctanh{\mathrm{arctanh}}

\def\O{{\cal O}}

\newcounter{theremark}
\newcommand{\rem}{\par\noindent\refstepcounter{theremark}\textbf{Remark \arabic{theremark}.} }
\newcommand\coef[3]{\beta^{#1}_{\!#2}[#3]}
\title{Series expansions for the third incomplete elliptic integral via partial fraction decompositions}
\author{D. Karp\footnote{Institute of Applied Mathematics, Vladivostok, Russia, e-mail:\emph{dmkrp@yandex.ru}},
A. Savenkova\footnote{Far Eastern State University, Vladivostok,
Russia, e-mail:\emph{asya-savenkova@yandex.ru}},
S.~M.~Sitnik\footnote{Voronezh Institute of the Ministry of
Internal Affairs of the Russian Federation,
e-mail:\emph{mathsms@yandex.ru}}}
\date{}
\begin{document}
\maketitle

\begin{center}
\parbox{12cm}{
\small\textbf{Abstract.}  We find convergent double series
expansions for Legendre's third incomplete elliptic integral valid
in overlapping subdomains of the unit square. Truncated expansions
provide asymptotic approximations in the neighbourhood of the
logarithmic singularity $(1,1)$ if one of the variables approaches
this point faster than the other.  Each approximation is
accompanied by an error bound.  For a curve with an arbitrary
slope at $(1,1)$ our expansions can be rearranged into asymptotic
expansions depending on a point on the curve.  For reader's
convenience we give some numeric examples and explicit expressions
for low-order approximations.}
\end{center}

\bigskip


\bigskip


\bigskip

\paragraph{1. Introduction.} Legendre's incomplete elliptic integral (EI) of the third kind is
defined by \cite[formula 13.6(3)]{Bat3}
\begin{equation}\label{eq:Pi-defined}
\Pi(\lambda,\nu,k)=\int\limits_{0}^{\lambda}\frac{dt}{(1+\nu{t^2})\sqrt{(1-t^2)(1-k^2t^2)}}.
\end{equation}
It is one of the three canonical forms given by Legendre in terms
of which all elliptic integrals can be expressed. We will only
consider the most important case $0\leq{k}\leq{1}$,
$0\leq{\lambda}\leq{1}$ and $\nu>-1$. Two other canonical forms
can be expressed in terms of $\Pi(\lambda,\nu,k)$: setting $\nu=0$
gives the first incomplete elliptic integral, while setting
$\nu=-k^2$ gives the second plus an elementary function (see
\cite[formulas 13.6(11)-13.6(13)]{Bat3}). When $\lambda=1$ the rhs
of (\ref{eq:Pi-defined}) is called the third complete elliptic
integral.

Series expansions, asymptotic approximations and inequalities for
the third incomplete elliptic integral have been studied by many
authors. Major contributions were made by Radon \cite{Radon},
Carlson \cite{Carlson1,CarlsonBook}, Carlson and Gustafson
\cite{Carlson2} and Lopez \cite{Lopez,Lopez1}. Carlson showed in
\cite{Carlson1} that $\Pi(\lambda,\nu,k)$ can be expressed in
terms of Lauricella hypergeometric function $F_D$ of three
variables (see \cite{Bat1}). In the same paper he noted that one
can derive rapidly convergent expansions by first expressing
Legendre's incomplete EIs in a different form. This form later
became known as symmetric standard EIs.   Today, most authors
consider asymptotic approximations for these symmetric EIs. Since
Legedre's incomplete elliptic integrals are connected with
symmetric EIs by certain simple relations (see
\cite{Carlson1,CarlsonBook}), it is possible to reformulate
expansions presented here in terms of symmetric elliptic
integrals.

In this paper we give convergent double series expansion for the
third incomplete elliptic integral (\ref{eq:Pi-defined}) valid in
two overlapping subregions of the unit square $[0,1]\times[0,1]$
in $(k^2,\lambda^2)$-plane defined by
\begin{equation}\label{eq:k-lambda-ineq}
\text{Region I}=\{k^2,\lambda^2\in[0,1]:
\lambda^2<1/(2-k^2)\Leftrightarrow
(1-k^2)\lambda^2/(1-\lambda^2)<1\},
\end{equation}
\begin{equation}\label{eq:region3}
\text{Region
II}=\{k^2,\lambda^2\in[0,1]:\lambda^2>2-1/k^2\Leftrightarrow
(1-\lambda^2)k^2/(1-k^2)<1\}.
\end{equation}
Combined these subregions cover the unit square completely.
Truncating our expansions one gets asymptotic approximations for
$(1-\lambda)/(1-k)\to{0}$ or $(1-\lambda)/(1-k)\to{\infty}$, i.e.
for curves with endpoint on the top or right side of the unit
square or with endpoint $\lambda=k=1$ and horizontal or vertical
tangent at that point. In all cases we give explicit bounds for
the remainder term.

For a given smooth curve approaching the singular point $(1,1)$
from any direction our double series expansions can be rearranged
into asymptotic expansions, where a point on the curve plays the
role of an asymptotic parameter.  In the last section of the paper
we illustrate our results by numerical examples.  Rather good
precision (with less than one percent relative error)  is reached
in most cases already by the first asymptotic term.  For reader's
convenience we also give explicit expressions for first several
terms of each expansion.

\setcounter{section}{2}
\paragraph{2.~Partial fraction decompositions.}
For a non-negative integer $n$ and $\alpha=0,1$ define
\begin{equation}\label{eq:phi-psi}
\phi_n^{\alpha}(t)=(1-t)^{-n}+(-1)^{\alpha}(1+t)^{-n},~~~\psi_n^{\alpha}(t)=\left\{\begin{array}{ll}[(1-t)^{-n}-(-1)^{\alpha}(1+t)^{-n}-2\alpha]/n,&
n>0,\\(-1)^{\alpha}\ln(1+t)-\ln(1-t),& n=0.
\end{array}\right.
\end{equation}
The second definition is consistent in the sense that
$\lim\limits_{n\to{0}}\psi_n^{\alpha}(t)=(-1)^{\alpha}\ln(1+t)-\ln(1-t)$,
 $\alpha=0,1$. Our first lemma is straightforward.
\begin{lemma}\label{lemma:phipsi-int}
The functions $\phi^{\alpha}_n$ and $\psi^{\alpha}_n$ are related
by{\em:}
\begin{equation}\label{eq:phi-int}
\int\limits_{0}^{\lambda}\phi^{\alpha}_n(t)dt=\psi^{\alpha}_{n-1}(\lambda),~~~~~n=1,2,\ldots,~~~\alpha=0,1.
\end{equation}
\end{lemma}

\begin{lemma}\label{lemma:pfd}
Let $n<2j$ be non-negative integers. Then the following partial
fraction decomposition holds true\emph{:}
\begin{equation}\label{eq:mn-frac}
\frac{t^{n}}{(1-t^2)^{j}}=\coef{n}{j}{1}\phi_1^{\alpha}(t)+
\coef{n}{j}{2}\phi_2^{\alpha}(t)+\cdots+
\coef{n}{j}{j}\phi_j^{\alpha}(t),
\end{equation}
where $\alpha=n\!\!\!\mod{\!2}$, $\coef{n}{j}{j}=2^{-j}$ and the
numbers $\coef{n}{j}{j-k}$, $k=1,2,\ldots,j-1$, are found from the
following recurrence \emph{(}here $\binom{n}{k}=0$ when
$k>n$\emph{):}
\begin{equation}\label{eq:alpha-recurrence}
\coef{n}{j}{j-k}=\frac{(-1)^k}{2^{j}}\binom{n}{k}-\sum\limits_{i=0}^{k-1}
(-2)^{i-k}\binom{j}{k-i}\coef{n}{j}{j-i},~~~~k=1,2,\ldots,j-1,
\end{equation}
or explicitly from
\begin{equation}\label{eq:beta-explicit}
\coef{n}{j}{j-k}=2^{-j-k}\sum_{i=0}^{k}(-2)^i\binom{n}{i}\binom{k+j-1-i}{j-1},
\quad 0\leq k\leq j.
\end{equation}
\end{lemma}

\textbf{Proof.} For computing the coefficients of the partial
fraction decomposition we use the trick routinely applied for
calculating residues.  Writing
\[
\frac{t^{n}}{(1-t)^{j}(1+t)^j}=\frac{\coef{n}{j}{1}}{1-t}+
\frac{\coef{n}{j}{2}}{(1-t)^2}+\cdots+
\frac{\coef{n}{j}{j}}{(1-t)^{j}}+(-1)^n\frac{R(t)}{(1+t)^j},
\]
where
\[
\frac{R(t)}{(1+t)^j}=\frac{\coef{n}{j}{1}}{1+t}+\frac{\coef{n}{j}{2}}{(1+t)^2}+\cdots+\frac{\coef{n}{j}{j}}{(1+t)^{j}},
\]
and multiplying throughout by $(1-t^2)^{j}$, we get:
\[
t^{n}=(1+t)^{j}(\coef{n}{j}{j}+\coef{n}{j}{j-1}(1-t)+\cdots+\coef{n}{j}{1}(1-t)^{j-1})+(1-t)^{j}R(t).
\]
Setting $t=1$ immediately reveals $\coef{n}{j}{j}=2^{-j}$.
Differentiate $k$ times and set $t=1$:
\[
[t^{n}]^{(k)}_{\vert t=1}=\sum\limits_{i=0}^{k}\binom{k}{i}
[(1+t)^{j}]^{(k-i)}_{\vert
t=1}\coef{n}{j}{j-i}(-1)^ii!,~~~~k=1,2,\ldots,j-1.
\]
We have:
\[
[t^{n}]^{(k)}_{\vert
t=1}=\left\{\begin{array}{lr}n!/(n-k)!,&~~k\leq{n}
\\
0,&~~k>n\end{array}\right.,~~~[(1+t)^{j}]^{(k-i)}_{\vert
t=1}=\frac{j!}{(j-k+i)!}2^{j-k+i}.
\]
Thus the numbers $\coef{n}{j}{j-k}$ can be found via the following
recurrence  (starting with $\coef{n}{j}{j}=2^{-j}$ and counting
downwards):
\[
\coef{n}{j}{j-k}(-1)^kk!2^{j}=\frac{n!}{(n-k)!}-\sum\limits_{i=0}^{k-1}\binom{k}{i}
\frac{j!}{(j-k+i)!}2^{j-k+i}\coef{n}{j}{j-i}(-1)^ii!,~~~~k=1,2,\ldots,j-1.
\]
A simple rearrangement of this formula yields
(\ref{eq:alpha-recurrence}) which can also be written as
$$
\sum_{i=0}^{k}(-2)^{i-k}\binom{j}{k-i}\coef{n}{j}{j-i}=\frac{(-1)^k}{2^j}\binom{n}{k}.
$$
This recurrence can be solved for $\coef{n}{j}{j-i}$ using the
inversion formula (see \cite{Riordan})
$$
a_n=\sum_{k=0}^{n}(-1)^{k}\binom{m}{n-k}d_k \quad \Rightarrow
\quad d_n = \sum_{k=0}^{n}(-1)^{k}\binom{n+m-1-k}{n-k} a_k,
$$
which leads to (\ref{eq:beta-explicit}).~$\square$

In the sequel we will only need even  $n$ therefore we put
$\phi_n(t)=\phi_n^0(t)$ and $\psi_n(t)=\psi_n^0(t)$ to simplify
notation.

\begin{lemma} The following partial fraction decomposition holds true for
$j=0,1,\ldots${\em:}
\begin{equation}\label{eq:mainpfd}
\frac{t^{2j}}{(1+\nu{t^2})(1-t^2)^{j+1}}=\frac{(-1)^j\nu}{(1+\nu)^{j+1}(1+{\nu}t^2)}+
\sum\limits_{k=0}^{j}\frac{(-1)^{k}\phi_{j+1-k}(t)}{(1+\nu)^{k+1}}\sum\limits_{i=0}^{k}a_{j+1-k,i}(1+\nu)^i,
\end{equation}
where the numbers $a_{n,m}$ are expressed in terms of the numbers
$\coef{k}{p}{i}$ defined by {\em(\ref{eq:alpha-recurrence})} or
{\em(\ref{eq:beta-explicit})} by
\begin{equation}\label{eq:a-defined}
a_{n,m}=(-1)^{m}\coef{2n+2m-4}{n+m}{n},~~n=1,\ldots,j+1;
~~m=\left\{\!\!\begin{array}{ll}1,\ldots,j, & n=1\\
0,\ldots,j+1-n,& n>1
\end{array}\right.\!\!; ~~a_{1,0}=\frac{1}{2}.
\end{equation}
\end{lemma}
\rem This lemma contains two statements: one is formula
(\ref{eq:a-defined}) for the coefficients of expansion
(\ref{eq:mainpfd}), and the other is that the numbers $a_{n,m}$
{\em are independent of j}. That is when $j$ increases we do not
need to update all the numbers $a_{n,m}$, instead we keep all
previously calculated numbers and complement them with new $j+1$
numbers with indices summing up to $j+1$: $n+m=j+1$.

\textbf{Proof}.  We will use induction in $j$. For $j=0$ we can
verify directly that:
\[
\frac{1}{(1+\nu{t^2})(1-t^2)}=\frac{\nu}{(1+\nu)(1+{\nu}t^2)}+\frac{1}{2(1+\nu)}
\left[\frac{1}{1-t}+\frac{1}{1+t}\right].
\]
Suppose now that (\ref{eq:mainpfd}) holds for a fixed $j$.  Then
for $j+1$,
\[
\frac{t^{2j+2}}{(1+\nu{t^2})(1-t^2)^{j+2}}=\frac{t^{2j}}{(1+\nu)(1-t^2)^{j+2}}-\frac{t^{2j}}{(1+\nu)(1+\nu{t^2})(1-t^2)^{j+1}}.
\]
Now using (\ref{eq:mn-frac}) with $n=2j$, $m=j+2$ for the first
term and (\ref{eq:mainpfd}) for the second gives after collecting
coefficients at $\phi_i(t)$:
\begin{eqnarray*}
&&\frac{t^{2j+2}}{(1+\nu{t^2})(1-t^2)^{j+2}}=\frac{(-1)^{j+1}\nu}{(1+\nu)^{j+2}(1+{\nu}t^2)}+\nonumber\\
&&+\frac{(-1)^{j+1}\left(a_{1,0}+a_{1,1}(1+\nu)+\cdots+a_{1,j}(1+\nu)^j\right)+\coef{2j}{j+2}{1}(1+\nu)^{j+1}}{(1+\nu)^{j+2}}
\phi_1(t)+\\
&&+\frac{(-1)^{j}\left(a_{2,0}+a_{2,1}(1+\nu)+\cdots+a_{2,j-1}(1+\nu)^{j-1}\right)+\coef{2j}{j+2}{2}(1+\nu)^{j-1}}{(1+\nu)^{j}}
\phi_2(t)+\\
&&+\cdots+\frac{(-1)a_{j+1,0}+\coef{2j}{j+2}{j+1}(1+\nu)}{(1+\nu)^2}\phi_{j+1}(t)
+\frac{\coef{2j}{j+2}{j+2}}{1+\nu}\phi_{j+2}(t).
\end{eqnarray*}
Now define
\[
a_{1,j+1}=(-1)^{j+1}\coef{2j}{j+2}{1},~~a_{2,j}=(-1)^{j}\coef{2j}{j+2}{2},~\ldots,~a_{j+2,0}=\coef{2j}{j+2}{j+2}.
\]
This shows that expansion for $j+1$ has the form
(\ref{eq:mainpfd}) with numbers $a_{n,m}$ given by
(\ref{eq:a-defined}). $\square$

Formula (\ref{eq:mainpfd}) combined with
Lemma~\ref{lemma:phipsi-int} leads to the following evaluation
\begin{equation}\label{eq:mainfracint}
\int\limits_{0}^{\lambda}\frac{t^{2j}dt}{(1+\nu{t^2})(1-t^2)^{j+1}}=\frac{(-1)^j\sqrt{\nu}}{(1+\nu)^{j+1}}\arctan(\lambda\sqrt{\nu})+
\sum\limits_{n=0}^{j}d_{j,n+1}(\nu)\psi_{n}(\lambda),
\end{equation}
where we introduced the notation
\begin{equation}\label{eq:d-defined}
d_{j,n+1}(\nu)=\frac{(-1)^{j-n}}{(1+\nu)^{j-n+1}}\sum\limits_{i=0}^{j-n}a_{n+1,i}(1+\nu)^i,~~n=0,1,2,\ldots,j.
\end{equation}
Here and henceforth $\sqrt{\nu}\arctan(\lambda\sqrt{\nu})$ will be
understood as $(\sqrt{-\nu})\arctanh(\lambda\sqrt{-\nu})$ for
negative $\nu$.

Another way of decomposing the left hand side of
(\ref{eq:mainpfd}) is given in the next lemma.
\begin{lemma}\label{lm:main-new}
The following expansions hold true for $j=0,1,\ldots${\em:}
\begin{equation}\label{eq:main-new}
\frac{t^{2j}}{(1+\nu{t^2})(1-t^2)^{j+1}}=\left[\frac{\nu}{1+\nu}\right]^{j+1}\frac{t^{2j}}{1+\nu{t^2}}
+\sum\limits_{n=1}^{j+1}\frac{\nu^{j-n+1}}{(1+\nu)^{j-n+2}}\frac{t^{2j}}{(1-t^2)^n},
\end{equation}
and
\begin{equation}\label{eq:forPi-new}
\frac{t^{2j}}{1+\nu{t^2}}=\frac{(-1/\nu)^j}{1+\nu{t^2}}-\sum\limits_{i=0}^{j-1}(-1/\nu)^{j-i}t^{2i}.
\end{equation}
\end{lemma}

\textbf{Proof}.  Both expansions can be verified by using the
summation formula for a finite geometric progression.~$\square$

Euler's integral representation for the Gauss hypergeometric
function ${_2F_1}$ (after simple variable change) combined with
(\ref{eq:phi-int}) and (\ref{eq:mn-frac}) gives
\begin{equation}\label{eq:ourEuler}
\int\limits_{0}^{\lambda}\frac{t^{2j}dt}{(1-t^2)^{n}}=\frac{\lambda^{2j+1}}{2j+1}
{_{2}F_{1}}(n,j+1/2;j+3/2;\lambda^2)=\sum\limits_{i=1}^{n}\coef{2j}{n}{i}\psi_{i-1}(\lambda).
\end{equation}
Lemma~\ref{lm:main-new} and formula (\ref{eq:ourEuler}) provide
two alternative ways to evaluate the integral in
(\ref{eq:mainfracint}).  We have
\begin{equation}\label{eq:mainfracint-new0}
\int\limits_{0}^{\lambda}\frac{t^{2j}dt}{(1+\nu{t^2})(1-t^2)^{j+1}}
=\left[\frac{\nu}{1+\nu}\right]^{j+1}\int\limits_{0}^{\lambda}\frac{t^{2j}dt}{1+\nu{t^2}}+
\frac{\lambda^{2j+1}}{2j+1}\sum\limits_{n=1}^{j+1}\frac{\nu^{j-n+1}{_{2}F_{1}}(n,j+1/2;j+3/2;\lambda^2)}{(1+\nu)^{j-n+2}}
\end{equation}
\begin{equation}\label{eq:mainfracint-new}
=\frac{(-1)^j\sqrt{\nu}}{(1+\nu)^{j+1}}\arctan(\lambda\sqrt{\nu})
+\frac{(-1)^j}{(1+\nu)^{j+1}}\sum\limits_{i=0}^{j-1}(-\nu)^{i+1}\frac{\lambda^{2i+1}}{2i+1}
+\sum\limits_{n=1}^{j+1}\frac{\nu^{j-n+1}}{(1+\nu)^{j-n+2}}\sum\limits_{i=1}^{n}\coef{2j}{n}{i}\psi_{i-1}(\lambda).
\end{equation}


\paragraph{3. Expansions for the third incomplete elliptic integral.}
The following relation will be of great help:
\begin{equation}\label{eq:Pi-relation}
\Pi(\lambda,\nu,k)=\Pi(\nu,k)-\frac{1}{(1+\nu)\sqrt{1-k^2}}\Pi\left(-\nu/(1+\nu),\sqrt{1-\lambda^2},\sqrt{-k^2/(1-k^2)}\right),
\end{equation}
where $\Pi(\nu,k)$ is the complete elliptic integral of the third
kind. Here we choose the branch of the first square root that is
positive for positive values of $1-\lambda^2$. The choice of the
branch of the second square root is immaterial since $\Pi$ depends
on the squared second argument only.  This relation can be easily
verified by representing the integral over $(0,\lambda)$ from
(\ref{eq:Pi-defined}) as the difference of integrals over $(0,1)$
and $(\lambda,1)$ and introducing the new integration variable
$u^2=1-t^2$.

\begin{theo}\label{th:Pi-first}
Let $k,\lambda\!\in\!\emph{Region I}$ and $\nu>-1$. Then for any
positive integer $N$ the following expansions hold true
\begin{equation}\label{eq:Pi-series}
\Pi(\lambda,\nu,k)=\!\sum\limits_{j=0}^{N-1}(-1)^j\frac{(1/2)_j}{j!}(1-k^2)^{j}
\left[\frac{(-1)^j\sqrt{\nu}}{(1+\nu)^{j+1}}\arctan(\lambda\sqrt{\nu})+
\sum\limits_{n=0}^{j}d_{j,n+1}(\nu)\psi_{n}(\lambda)\right]\!+\!R_{1,N}(k,\lambda,\nu),
\end{equation}
and
\[
\Pi(\lambda,\nu,k)=\sum\limits_{j=0}^{N-1}\frac{(1/2)_j}{j!}(1-k^2)^j
\left\{\frac{\sqrt{\nu}}{(1+\nu)^{j+1}}\arctan(\lambda\sqrt{\nu})\right.
\]
\begin{equation}\label{eq:Pi-new}
\left.+\frac{1}{(1+\nu)^{j+1}}\sum\limits_{i=0}^{j-1}(-\nu)^{i+1}\frac{\lambda^{2i+1}}{2i+1}
+\sum\limits_{n=1}^{j+1}\frac{(-1)^j\nu^{j-n+1}}{(1+\nu)^{j-n+2}}\sum\limits_{i=1}^{n}\coef{2j}{n}{i}\psi_{i-1}(\lambda)\right\}
\!+\!R_{1,N}(k,\lambda,\nu),
\end{equation}
where the functions $\psi_{n}(\lambda)$  and $d_{j,n}(\nu)$ are
defined by {\em(\ref{eq:phi-psi})} and {\em(\ref{eq:d-defined})},
respectively, and the remainder term satisfies
\begin{equation}\label{eq:R1Piest1}
|R_{1,N}(k,\lambda,\nu)|\leq
\frac{(1/2)_N\lambda}{2\min(1,1+\nu)NN!}\left[\frac{(1-k^2)\lambda^2}{1-\lambda^2}\right]^{N}.
\end{equation}
\end{theo}

\textbf{Proof.} Put ${k^\prime}^{2}=1-k^2$ and calculate using
binomial expansion and termwise integration:
\[
\Pi(\lambda,\nu,k)=\int\limits_{0}^{\lambda}\frac{dt}{(1+\nu{t^2})\sqrt{(1-t^2)(1-k^2t^2)}}=
\int\limits_{0}^{\lambda}\!\!\!\frac{dt}{(1+\nu{t^2})(1-t^2)}\left(\!1+\frac{{k^\prime}^2t^2}{1-t^2}\right)^{-1/2}
\]\[
=\int\limits_{0}^{\lambda}\!\!\!\frac{dt}{(1+\nu{t^2})(1-t^2)}
\left(\sum\limits_{j=0}^{N-1}(-1)^j\frac{(1/2)_j}{j!}\frac{{k^\prime}^{2j}t^{2j}}{(1-t^2)^j}+
\sum\limits_{j=N}^{\infty}(-1)^j\frac{(1/2)_j}{j!}\frac{{k^\prime}^{2j}t^{2j}}{(1-t^2)^j}\right)
\]
\begin{equation}\label{eq:Pi-exp-in-proof}
=\sum\limits_{j=0}^{N-1}(-1)^j\frac{(1/2)_j}{j!}{k^\prime}^{2j}
\int\limits_{0}^{\lambda}\!\!\!\frac{t^{2j}dt}{(1+\nu{t^2})(1-t^2)^{j+1}}+R_{1,N}(k,\lambda,\nu),
\end{equation}
where the remainder term is given by
\begin{equation}\label{eq:R1}
R_{1,N}(k,\lambda,\nu)=\sum\limits_{j=N}^{\infty}(-1)^j\frac{(1/2)_j}{j!}(1-k^2)^{j}
\int\limits_{0}^{\lambda}\!\!\!\frac{t^{2j}dt}{(1+\nu{t^2})(1-t^2)^{j+1}}.
\end{equation}
Using representations (\ref{eq:mainfracint}) and
(\ref{eq:mainfracint-new}) for the integral in
(\ref{eq:Pi-exp-in-proof}) we get (\ref{eq:Pi-series}) and
(\ref{eq:Pi-new}), respectively. To prove the error bound
(\ref{eq:R1Piest1}) we note that the series (\ref{eq:R1}) has
alternating signs. The following estimate shows that the terms in
(\ref{eq:R1}) monotonically decrease in absolute value:
\[
\frac{(1/2)_{j+1}}{(j+1)!}(1-k^2)^{j+1}\int\limits_{0}^{\lambda}\frac{t^{2j+2}dt}{(1+\nu{t^2})(1-t^2)^{j+2}}=
\frac{(1/2)_j(1/2+j)}{j!(j+1)}\int\limits_{0}^{\lambda}\frac{(1-k^2)^{j}t^{2j}}{(1+\nu{t^2})(1-t^2)^{j+1}}\frac{t^2(1-k^2)}{(1-t^2)}dt
\]\[
\leq\frac{(1/2)_j}{j!}(1-k^2)^{j}\int\limits_{0}^{\lambda}\frac{t^{2j}}{(1+\nu{t^2})(1-t^2)^{j+1}}\frac{\lambda^2(1-k^2)}{(1-\lambda^2)}dt\leq
\frac{(1/2)_j}{j!}(1-k^2)^{j}\int\limits_{0}^{\lambda}\frac{t^{2j}dt}{(1+\nu{t^2})(1-t^2)^{j+1}}.
\]
The last inequality is due to (\ref{eq:k-lambda-ineq}).   Hence we
are in the position to apply the Leibnitz convergence test and the
remainder term does not exceed the first term in (\ref{eq:R1}):
\[
|R_{1,N}(k,\lambda,\nu)|\leq\frac{(1/2)_N}{N!}(1-k^2)^{N}\int\limits_{0}^{\lambda}\frac{t^{2N}dt}{(1+\nu{t^2})(1-t^2)^{N+1}}
\leq
\frac{(1/2)_N(1-k^2)^{N}}{N!\min(1,1+\nu)}\int\limits_{0}^{\lambda}\frac{t^{2N}dt}{(1-t^2)^{N+1}}.
\]
The integral on the right hand side satisfies the following
asymptotically precise  (as $\lambda\to{1}$) estimate
\begin{equation}\label{eq:lambdaint_est}
f_1(\lambda)\equiv\int\limits_{0}^{\lambda}\frac{t^{2N}dt}{(1-t^2)^{N+1}}\leq\frac{\lambda^{2N+1}}{2N(1-\lambda^2)^N}\equiv{f_2(\lambda)},
\end{equation} which is valid for all $\lambda\in(0,1)$ and $N>0$
(not necessarily integer).  Indeed, $f_1(0)=f_2(0)=0$ and
\[
\frac{f_1'(\lambda)}{f_2'(\lambda)}=\frac{2N}{2N+1-\lambda^2}<1,
~~~\lambda\in(0,1).
\]
This estimate immediately leads to (\ref{eq:R1Piest1}).~$\square$

\rem\label{rm:ordering}  The error bound (\ref{eq:R1Piest1}) shows
that expansions (\ref{eq:Pi-series}) and (\ref{eq:Pi-new}) are
asymptotic for $\lambda\to{1}$ with either constant $k$ or
$k\to{1}$ in a way that $(1-\lambda)/(1-k)\to{0}$, so that we can
approach the singular point $(1,1)$ along a curve in $(k,\lambda)$
plane having infinite slope at $(1,1)$. These expansions can be
transformed into asymptotic expansion for $k,\lambda\to{1}$ along
an arbitrary curve if an ordering can be introduced into the
matrix of functions forming the expansion that converts this
matrix into an asymptotic scale. Indeed, we have the matrix
\vspace{-0.2in}
\[
\begin{array}{cccccc}
&&&&&\\
\psi_0(\lambda) &  &  & &\\
\downarrow&&&&&\\
\arctan(\lambda\sqrt{\nu}) &  &  & &\\
\downarrow&&&&&\\
(1-k^2)\arctan(\lambda\sqrt{\nu})& (1-k^2)\psi_0(\lambda) & \leftarrow(1-k^2)\psi_1(\lambda) & &  &   \\
\downarrow&\downarrow&\downarrow&&\\
(1-k^2)^2\arctan(\lambda\sqrt{\nu})& (1-k^2)^2\psi_0(\lambda) & \leftarrow(1-k^2)^2\psi_1(\lambda) & \leftarrow(1-k^2)^2\psi_2(\lambda) &  \\
\downarrow&\downarrow&\downarrow&\downarrow&\\
\cdots & \cdots & \cdots & \cdots \\
\end{array}
\]
Arrows indicate the ordering that is true no matter from which
direction we approach the point $(1,1)$.  For instance,
$(1-k^2)^2\psi_1(\lambda)=o((1-k^2)\psi_1(\lambda))$ as
$k,\lambda\to{1}$. However, unless we specify a curve in the
$(k,\lambda)$ plane this is only a partial ordering. For instance,
we do not know how to compare $(1-k^2)^2\psi_1(\lambda)$ and
$(1-k^2)^3\psi_2(\lambda)$. Choosing a curve possessing a tangent
at the point $(1,1)$ introduces the complete ordering into the
above matrix. Rearranging terms in the expansions
(\ref{eq:Pi-series}) or (\ref{eq:Pi-new}) according to this
ordering converts it into the asymptotic expansion for the
specified curve. If the slope of our curve at $(1,1)$ is strictly
between $0$ and $\pi/2$, then all the diagonals of the above
matrix will be of the same asymptotic order and the order will be
decreasing as we move away from the main diagonal.

\rem The error bound (\ref{eq:R1Piest1}) shows that expansions
(\ref{eq:Pi-series}) and (\ref{eq:Pi-new}) are convergent for $k$,
$\lambda$ in $\text{Region I}$. Taking $N=\infty$ in
(\ref{eq:Pi-series}) and (\ref{eq:Pi-new}) and summing up we
obtain convergent double series representations
\begin{equation}\label{eq:Pi-series1}
\Pi(\lambda,\nu,k)=\sqrt{\frac{\nu}{1+\nu}}(\nu+k^2)^{-1/2}\arctan(\lambda\sqrt{\nu})
+
\sum\limits_{j=0}^{\infty}\sum\limits_{n=0}^{j}\frac{(1/2)_j(1-k^2)^{j}}{(-1)^jj!}
d_{j,n+1}(\nu)\psi_{n}(\lambda),
\end{equation}
\begin{equation}\label{eq:Pi-new1}
\Pi(\lambda,\nu,k)=\int\limits_{0}^{\lambda}\frac{[\nu/(1+\nu)]dt}{(1+\nu{t^2})\sqrt{1+\mu{t^2}}}
+\sum\limits_{j=0}^{\infty}\sum\limits_{n=1}^{j+1}
\frac{(1/2)_j\lambda^{2j+1}(1-k^2)^j\nu^{j-n+1}}{(-1)^jj!(2j+1)(1+\nu)^{j-n+2}}
{_{2}F_{1}}(n,j+1/2;j+3/2;\lambda^2),
\end{equation}
where $\mu=\nu(1-k^2)/(1+\nu)$. The first series requires
additional condition $-k^2<\nu$ for convergence.


To derive an expansion valid in $\text{Region II}$ we start with
an expansion valid in a neighbourhood of $(0,0)$ given by Carlson
\cite{Carlson1}:
\begin{equation}\label{eq:Pi-Carlson}
\Pi(\lambda,\nu,k)=\sum\limits_{m=0}^{\infty}\frac{\lambda^{2m+1}}{2m+1}\sum\limits_{n=0}^{m}(-\nu)^{m-n}\frac{(1/2)_n}{n!}{_2F_1}(-n,1/2;1/2-n;k^2).
\end{equation}
Using
\begin{equation}\label{eq:2F1}
\frac{(1/2)_n}{n!}{_2F_1}(-n,1/2;1/2-n;k^2)={_2F_1}(-n,1/2;1;1-k^2),
\end{equation}
which is a limiting case of the well-known analytic extension
formula for ${_2F_1}$ (see \cite[formula 2.10(1)]{Bat1}), we
obtain the expansion
\begin{equation}\label{eq:Pi-Gauss}
\Pi(\lambda,\nu,k)=\sum\limits_{m=0}^{\infty}\frac{\lambda^{2m+1}}{2m+1}\sum\limits_{n=0}^{m}(-\nu)^{m-n}{_2F_1}(-n,1/2;1;1-k^2),
\end{equation}
which is valid for $|\nu|\lambda^2<1$ and $|k|<1/\lambda$ as shown
by Carlson. An application of (\ref{eq:Pi-relation}) leads to

\begin{theo}
For $k,\lambda\!\in\!\emph{Region II}$, $\nu$ in the range
\begin{equation}\label{eq:nu-lambda-ineq}
((1-\lambda^2)|\nu|)/(\nu+1)<1,
\end{equation}
and positive integer $N$ the following expansion holds true
\begin{equation}\label{eq:Pi-Gauss-reverse}
\Pi(\lambda,\nu,k)\!=\!\Pi(\nu,k)
-\sqrt{\frac{1-\lambda^2}{1-k^2}}
\sum\limits_{m=0}^{N-1}\frac{(1-\lambda^2)^{m}}{2m+1}\sum\limits_{n=0}^{m}\frac{\nu^{m-n}{_2F_1}(-n,1/2;1;(1-k^2)^{-1})}{(1+\nu)^{m-n+1}}
+R_{2,N}(k,\nu,\lambda),
\end{equation}
with the error bound given by
\begin{equation}\label{eq:R2-bound}
|R_{2,N}(k,\nu,\lambda)|\leq\frac{(1-\lambda^2)^{N+1/2}}{2N+1}
\left[\max\left(\frac{k^2}{1-k^2},\frac{|\nu|}{\nu+1},1\right)\right]^{N+1/2}f(\lambda,k,\nu),
\end{equation}
where $f(\lambda,k,\nu)$ is defined by
\begin{equation}\label{eq:f}
f(\lambda,k,\nu)=\left\{
\begin{array}{ll}
\frac{1}{(1+\nu)k}\left[1-\frac{(1-\lambda^2)k^2}{1-k^2}\right]^{-1}\left[1-\frac{|\nu|(1-k^2)}{(\nu+1)k^2}\right]^{-1}\!\!\!\!
,~~\max\!\left(\!\frac{k^2}{1-k^2},\frac{|\nu|}{\nu+1}\right)=\frac{k^2}{1-k^2}>1;\\
\frac{1}{\sqrt{(1-k^2)(1+\nu)|\nu|}}\left[1-\frac{(1-\lambda^2)|\nu|}{1+\nu}\right]^{-1}
\left[1-\frac{(\nu+1)k^2}{|\nu|(1-k^2)}\right]^{-1}\!\!\!\!,
~~\frac{|\nu|}{\nu+1}>\frac{k^2}{1-k^2}>1;\\
\frac{1}{\sqrt{(1-k^2)(1+\nu)|\nu|}}
\left[1-\frac{(1-\lambda^2)|\nu|}{1+\nu}\right]^{-1}\left(\left[1-\frac{1+\nu}{|\nu|}\right]^{-1}+\frac{1}{\lambda^2}\right),
~~\frac{|\nu|}{\nu+1}>1\geq\frac{k^2}{1-k^2};\\
\frac{1}{(1+\nu)\sqrt{1-k^2}}\left[1-\frac{(1-\lambda^2)|\nu|}{1+\nu}\right]^{-1}
\left(\left[1-\frac{|\nu|}{1+\nu}\right]^{-1}+\frac{1}{\lambda^2}\right),
~~\max\!\left(\!\frac{k^2}{1-k^2},\frac{|\nu|}{\nu+1},1\right)=1;\\
\frac{1}{\sqrt{(1-k^2)(1+\nu)|\nu|}}\left[1-\frac{(1-\lambda^2)|\nu|}{(1+\nu)}\right]^{-1}
\left(\left[1-\frac{(1-\lambda^2)|\nu|}{(1+\nu)}\right]^{-1}+N\right),
~~\max\!\left(\!\frac{k^2}{1-k^2},1\right)=\frac{|\nu|}{\nu+1}.
\end{array} \right.
\end{equation}
\end{theo}

\textbf{Proof.} To derive the form of (\ref{eq:Pi-Gauss-reverse})
expand the second term on the righthand side of
(\ref{eq:Pi-relation}) into series (\ref{eq:Pi-Gauss}). Carlson's
conditions $|\nu|\lambda^2<1$ and $|k|<1/\lambda$ written for the
parameters of this second term become (\ref{eq:region3}) and
(\ref{eq:nu-lambda-ineq}). Now consider the remainder term
\[
R_N(k,\nu,\lambda)=\frac{1}{1+\nu}\left[\frac{1-\lambda^2}{1-k^2}\right]^{1/2}\tilde{R}_N(k,\nu,\lambda),
\]\[
\tilde{R}_N(k,\nu,\lambda)=\sum\limits_{m=N}^{\infty}\frac{(1-\lambda^2)^{m}}{2m+1}\sum\limits_{n=0}^{m}\left(\frac{\nu}{1+\nu}\right)^{m-n}
\!\!\!{_2F_1}(-n,1/2;1;(1-k^2)^{-1}).
\]
Change the order of summation to get
\begin{multline}\label{eq:RPi-form1}
\tilde{R}_N(k,\nu,\lambda)=\sum\limits_{n=0}^{\infty}\sum\limits_{s=0}^{\infty}\frac{(1-\lambda^2)^{s+\max(n,N)}}{2(s+\max(n,N))+1}\left(\frac{\nu}{1+\nu}\right)^{s+\max(n,N)-n}
\!\!\!{_2F_1}(-n,1/2;1;(1-k^2)^{-1})
\\
=\sum\limits_{n=0}^{\infty}(1-\lambda^2)^n{_2F_1}(-n,1/2;1;(1-k^2)^{-1})\left[\frac{(1-\lambda^2)\nu}{1+\nu}\right]^{\max(0,N-n)}
\sum\limits_{s=0}^{\infty}
\frac{(1-\lambda^2)^{s}(\nu/(1+\nu))^s}{2(s+\max(n,N))+1}.
\end{multline}
Since ($L>-1/2$)
\[
\sum\limits_{s=0}^{\infty}\frac{x^s}{2s+2L+1}=\frac{x^{-L-1/2}}{2}\int\limits_{0}^{x}\frac{t^{L-1/2}}{1-t}dt,
\]
we may write
\begin{equation}\label{eq:RPi-form2}
\tilde{R}_N(k,\nu,\lambda)=\frac{1}{2}\left[\frac{1+\nu}{(1-\lambda^2)\nu}\right]^{1/2}
\!\!\!\!\int\limits_{0}^{\frac{(1-\lambda^2)\nu}{1+\nu}}\frac{dt}{\sqrt{t}(1-t)}
\sum\limits_{n=0}^{\infty}\left[\frac{1+\nu}{\nu}\right]^{n}t^{\max(n,N)}{_2F_1}(-n,1/2;1;(1-k^2)^{-1}).
\end{equation}
The series in (\ref{eq:RPi-form2}) can be broken into two parts:
\begin{equation}\label{eq:breakin2}
t^{N}\sum\limits_{n=0}^{N-1}\left[\frac{1+\nu}{\nu}\right]^{n}{_2F_1}(-n,1/2;1;(1-k^2)^{-1})+
\sum\limits_{n=N}^{\infty}\left[\frac{1+\nu}{\nu}\right]^{n}t^{n}{_2F_1}(-n,1/2;1;(1-k^2)^{-1}).
\end{equation}
To give a bound for (\ref{eq:RPi-form2}) we will need two
ingredients. First is the estimate
\begin{equation}\label{eq:2F1-est}
|{_2F_1}(-n,1/2;1;(1-k^2)^{-1})|\leq \left\{\!\!\begin{array}{ll}
|k^2/(1-k^2)|^n,~~1/2\leq{k^2}<1,\\1,~~0\leq{k^2}\leq{1/2},
\end{array}\right.
\end{equation}
which is obtained by writing ${_2F_1}$ as the Legendre polynomial
\cite[formula(7.3.176)]{Prud3}:
\begin{equation}\label{eq:2F1-Legendre}
{_{2}F_{1}}(-n,1/2;1;(1-k^2)^{-1})\!=\!\left[\frac{-k^2}{1-k^2}\right]^{\frac{n}{2}}\!\!P_n\left(\frac{1-2k^2}{2\sqrt{-k^2(1-k^2)}}\right)\!=
i^n\left[\frac{k^2}{1-k^2}\right]^{\frac{n}{2}}\!\!P_n\left(\frac{-i(1-2k^2)}{2k\sqrt{(1-k^2)}}\right)
\end{equation}
and applying the inequality (valid for all complex $z$)
$|P_n(z)|\leq |z+\sqrt{z^2-1}|^{n}$, where the branch of the
square root is chosen so that $|z+\sqrt{z^2-1}|\geq{1}$.  This
inequality follows immediately from the first Laplace integral for
$P_n(z)$ (see \cite[formula 3.7(6)]{Bat1}). The second ingredient
is the asymptotically precise (as $x\to{0}$) estimate (verified by
comparing derivatives of both sides and noticing that
$f(0)=g(0)$):
\begin{equation}\label{eq:int-est}
f(x)=\int\limits_{0}^x\frac{t^Mdt}{1-t}\leq
\frac{x^{M+1}}{(M+1)(1-x)}=g(x),
\end{equation}
valid for $M>-1$ and $0\leq{x}<1$. To prove
(\ref{eq:R2-bound})--(\ref{eq:f}) we need to consider the
following cases:

a) $1/2\leq{k^2}<1$, $(1+\nu)k^2/(|\nu|(1-k^2))>{1}$; b)
$1/2\leq{k^2}<1$, $(1+\nu)k^2/(|\nu|(1-k^2))<{1}$;

c) $1/2\leq{k^2}<1$, $(1+\nu)k^2/(|\nu|(1-k^2))={1}$; d)
$0\leq{k^2}\leq 1/2$, $\nu>-1/2$;

e) $0\leq{k^2}\leq 1/2$, $-1<\nu<-1/2$; and f) $0\leq{k^2}\leq
1/2$, $\nu=-1/2$.

\noindent Take a). Inequality (\ref{eq:2F1-est}) gives for the
first sum in (\ref{eq:breakin2}):
\[
|t|^{N}\sum\limits_{n=0}^{N-1}\left|\frac{1+\nu}{\nu}\right|^{n}|{_2F_1}(-n,1/2;1;(1-k^2)^{-1})|\leq
\left|\frac{t(1+\nu)k^2}{\nu(1-k^2)}\right|^{N}\frac{1}{[(\nu+1)k^2]/[|\nu|(1-k^2)]-1}.
\]
Similarly, the second sum in (\ref{eq:breakin2}) satisfies
\[
\left|\sum\limits_{n=N}^{\infty}\left[\frac{1+\nu}{\nu}\right]^{n}t^{n}{_2F_1}(-n,1/2;1;(1-k^2)^{-1})\right|
\leq
\frac{\left|\frac{t(1+\nu)k^2}{\nu(1-k^2)}\right|^{N}}{1-\left|\frac{t(1+\nu)k^2}{\nu(1-k^2)}\right|}
\leq\frac{\left|\frac{t(1+\nu)k^2}{\nu(1-k^2)}\right|^{N}}{1-\frac{(1-\lambda^2)k^2}{(1-k^2)}}.
\]
The last inequality is due to the fact that
$t\in[0,(1-\lambda^2)|\nu|/(\nu+1)]$. Hence the complete sum
(\ref{eq:breakin2}) is estimated by
\[
\left|\sum\limits_{n=0}^{\infty}\left[\frac{1+\nu}{\nu}\right]^{n}t^{\max(n,N+1)}{_2F_1}(-n,1/2;1;(1-k^2)^{-1})\right|
\leq
\left|\frac{t(1+\nu)k^2}{\nu(1-k^2)}\right|^{N}\frac{\frac{k^2(1-\lambda^2)}{(1-k^2)}\left[\frac{\nu+1}{|\nu|(1-\lambda^2)}-1\right]}
{\left[\frac{(\nu+1)k^2}{|\nu|(1-k^2)}-1\right]\left[1-\frac{(1-\lambda^2)k^2}{1-k^2}\right]}.
\]
Substituting this estimate into (\ref{eq:RPi-form2}) and using
(\ref{eq:int-est}) we get:
\[
|\tilde{R}_N(k,\nu,\lambda)|\leq
\frac{1+\nu}{(2N+1)|\nu|}\left[\frac{(1-\lambda^2)k^2}{(1-k^2)}\right]^{N}
\frac{k^2/(1-k^2)}{\left[\frac{(\nu+1)k^2}{|\nu|(1-k^2)}-1\right]\left[1-\frac{(1-\lambda^2)k^2}{1-k^2}\right]},
\]
which is equivalent to (\ref{eq:R2-bound}) plus first line of
(\ref{eq:f}). Proofs in all other cases require only minor
modifications.~~$\square$

\rem Error bound (\ref{eq:R2-bound}) shows that expansion
(\ref{eq:Pi-Gauss-reverse}) is convergent and asymptotic as
$k\to{1}$ with either constant $\lambda$ or $\lambda\to{1}$ in a
way that $(1-k)/(1-\lambda)\to{0}$, so that the point $(1,1)$ is
approached along a curve in $(k,\lambda)$ plane having zero slope
at $(1,1)$. It can be rearranged into an asymptotic expansion for
$k,\lambda\to{1}$ along an arbitrary curve as explained in
Remark~\ref{rm:ordering}.

\rem  Expansion (\ref{eq:Pi-Gauss-reverse}) can be derived without
any use of Carlson's expansion (\ref{eq:Pi-Carlson}) by expanding
the rightmost term in (\ref{eq:Pi-relation}) into binomial series,
similarly to the proof of Theorem~\ref{th:Pi-first}.

\rem The complete elliptic integral $\Pi(\nu,k)$ in
(\ref{eq:Pi-Gauss-reverse}) can be approximated using the
following asymptotic expansion valid for $k\to{1}$:
\begin{multline}\label{eq:Picomplete-asymp}
\Pi(\nu,k)\!=\sum\limits_{n=0}^{M-1}\frac{(1/2)_n(1/2)_n}{(n!)^2}(1-k^2)^n(\psi(1+n)-\psi(1/2+n)-\frac{1}{2}\ln(1-k^2)){_2F_1}(1,1/2+n;1/2;-\nu)\\
-\sum\limits_{n=0}^{M-1}\frac{(1/2)_n(1/2)_n}{2(n!)^2}(1-k^2)^n{_2F_1}^{(0,1,0,0)}(1,1/2+n;1/2;-\nu)+\O\left((1-k)^M\ln\frac{1}{1-k}\right).
\end{multline}
Here ${_2F_1}^{(0,1,0,0)}$ is the derivative of ${_2F_1}$ in the
second parameter.  This expansion can be obtained by writing
$\Pi(\nu,k)$ as Appell's hypergeometric function $F_1$ (see
\cite{Carlson1}), expressing $F_1$ as a series in terms of Gauss
hypergeometric functions ${_2F_1}$ and applying \cite[formula
2.10(12)]{Bat1} to ${_2F_1}$.

\setcounter{section}{4}
\paragraph{4. Results of computations.} In this section we give several examples of computations with
the expansions obtained above. Denote by
$\widetilde{\Pi}_N(\lambda,\nu, k)$ the $N$-th order approximation
in (\ref{eq:Pi-series}), i.e. the right hand side of
(\ref{eq:Pi-series}) without the remainder term
$R_{1,N}(\lambda,\nu, k)$. Then:
\begin{equation}
\widetilde{\Pi}_1(\lambda,\nu,k)=\frac{1}{2(1+\nu)}\ln\frac{1+\lambda}{1-\lambda}+
\frac{\sqrt{\nu}\arctan\left(\lambda\sqrt{\nu}\right)}{1+\nu},\nonumber
\end{equation}
\begin{gather}
\widetilde{\Pi}_2(\lambda,\nu,k)=\widetilde{\Pi}_1(\lambda,\nu,k)
-\frac{1-k^2}{8(1+\nu)(1-\lambda)}\left\{1+(1-\lambda)\frac{\nu-1}{\nu+1}\ln{\frac{1+\lambda}{1-\lambda}}
\right.\nonumber\\
\left.
-4\frac{1-\lambda}{1+\nu}\sqrt{\nu}\arctan\left(\lambda\sqrt{\nu}\right)
-\frac{1-\lambda}{1+\lambda}\right\},\nonumber
\end{gather}
\begin{gather}
\widetilde{\Pi}_3(\lambda,\nu, k) =
\widetilde{\Pi}_2(\lambda,\nu,k)+
{\frac{3}{128(1+\nu)}}\,\frac{(1-k^2)^2}{(1-\lambda)^2}
\left\{1-2\lambda\left(\frac{4}{1+\nu}+1\right)
\frac{1-\lambda}{1+\lambda}  \right.\nonumber\\
\qquad\left. +\ln{\frac {1+\lambda}{1-\lambda}} \left(
\frac{8}{(1+\nu)^2}- \frac{4}{1+\nu} - 1 \right)(1-\lambda)^2 +
16\,{\frac { \left( 1-\lambda \right) ^{2}\sqrt {\nu} \arctan
\left( \lambda\sqrt {\nu} \right) }{ \left( 1+\nu \right)
^{2}}}-{\frac { \left( 1-\lambda \right) ^{2}}{ \left( 1+\lambda
\right) ^{2}}} \right\},\nonumber
\end{gather}
\renewcommand\arraycolsep{3pt}
\setlength{\arrayrulewidth}{0.2pt}%
\setlength{\doublerulesep}{0pt}
{\small
$$
\begin{array}{!{\vrule width 1.2pt\relax}c|c|c|c!{\vrule width 1.2pt\relax}c|c|c!{\vrule width 1.2pt\relax}c|c|c!{\vrule width 1.2pt\relax}}
\hline\hline\hline\hline\hline &&&&&&&&&\\[-10pt]
\lambda & k & \displaystyle{\frac{1-k}{1-\lambda}} &
{\Pi}(\nu,\lambda, k) & \widetilde{\Pi}_1(\nu,\lambda,k) &
\begin{array}{c}\text{relative}\\[-3pt]\text{error}\end{array} &
\begin{array}{c}\text{relative}\\[-3pt]\text{error} \\[-3pt]\text{bound}\end{array} &
\widetilde{\Pi}_3(\lambda,\nu,k) &
\begin{array}{c}\text{relative}\\[-3pt]\text{error}\end{array} &
\begin{array}{c}\text{relative}\\[-3pt]\text{error} \\[-3pt]\text{bound}\end{array}
\\[5pt]
\hline\hline\hline\hline\hline
          .5 & .9 & .2 & .37138 & .37409 & -.730\!\!\times\!\!10^{-2}& .213\!\!\times\!\!10^{-1} &.37138 &  -.634\!\!\times\!\!10^{-5} & .178\!\!\times\!\!10^{-4}\\
\hline  .6 & .99 & .025 & .41973 & .42022 &-.117\!\!\times\!\!10^{-2} &  .400\!\!\times\!\!10^{-2}& .41973&  -.306\!\!\times\!\!10^{-7}& .104\!\!\times\!\!10^{-6}\\
\hline  .75 & .999 & .004 & .48662 & .48673 &-.218\!\!\times\!\!10^{-3} & .990\!\!\times\!\!10^{-3}& .48662& -.289\!\!\times\!\!10^{-9} & .136\!\!\times\!\!10^{-8}\\
\hline .9 & .99999 & .0001 & .57202 & .57203 &-.553\!\!\times\!\!10^{-5} & .335\!\!\times\!\!10^{-4}&.57202 & -.784\!\!\times\!\!10^{-14} & .508\!\!\times\!\!10^{-13}\\
\hline\hline\hline\hline\hline
\end{array}
$$
} \textbf{Table 4.1.} \textsl{Numerical examples for the
approximation \emph{(\ref{eq:Pi-series})} with $k\to{1}$, $\lambda
\to{1}$, $(1-k)/(1-\lambda)\to{0}$ and $\nu=7$. The sixth and the
ninth columns represents the relative errors
$R_{1,N}(\lambda,\nu,k)/\Pi(\lambda,\nu,k)$ in
\emph{(\ref{eq:Pi-series})}. The seventh and the tenths columns
represent relative error bounds as given by the rhs of
\emph{(\ref{eq:R1Piest1})} divided by $\Pi(\lambda,\nu,k)$.}

\medskip

For $k$ and $\lambda$ satisfying (\ref{eq:region3}) denote by
$\widehat{\Pi}_N(\lambda,\nu, k)$ the $N$-th order approximation
from (\ref{eq:Pi-Gauss-reverse}), i.e. the right hand side of
(\ref{eq:Pi-Gauss-reverse}) without remainder term
$R_{2,N}(\lambda,\nu, k)$. Then:
$$
\widehat{\Pi}_1(\lambda,\nu, k)= \Pi(\nu,k)
-\frac{1}{1+\nu}\left({\frac {1-\lambda^2}{1-k^2}}\right)^{1/2},
$$
$$
\widehat{\Pi}_2(\lambda,\nu, k)=\widehat{\Pi}_1(\lambda,\nu, k)
-\frac{1}{3}\, \left({\frac {1-\lambda^2}{1-k^2}}\right)^{3/2}
\left\{ \frac{1-k^2}{1+\nu} \left( 1+{ \frac {\nu}{1+\nu}} \right)
- \frac{1}{2}\, \frac{1}{1+\nu} \right\},
$$
\begin{multline}
\widehat{\Pi}_3(\lambda, \nu,k)=\widehat{\Pi}_2(\lambda,\nu, k) -
\frac{1}{5}\, \left( {\frac {1-{\lambda}^{2}}{1-{k}^{2}}} \right)
^{5/2} \left\{\frac{(1-k^2)^{2}}{1+\nu}\left(1+{\frac{\nu}{1+\nu}}+{\frac {{\nu}^{2}}{(1+\nu)^{2}}}\right)\right.\\
\left.-\frac{1}{2}\,\frac{1-k^2}{1+\nu}\left(2+{\frac{\nu}{1+\nu}}\right)+
\frac{3}{8}\,\frac{1}{1+\nu} \right\},\nonumber
\end{multline}
{\small
$$
\begin{array}{!{\vrule width 1.2pt\relax}c|c|c|c!{\vrule width 1.2pt\relax}c|c|c!{\vrule width 1.2pt\relax}c|c|c!{\vrule width 1.2pt\relax}}
\hline\hline\hline\hline\hline &&&&&&&&&\\[-10pt]
\lambda & k &
\displaystyle{\frac{1-\lambda}{1-k}}&\Pi(\lambda,\nu,k)&
\widehat{\Pi}_1(\lambda,\nu,k)&
\begin{array}{c}\text{relative}\\[-3pt]\text{error}\end{array}&
\begin{array}{c}\text{relative}\\[-3pt]\text{error} \\[-3pt]\text{bound}\end{array}
& \widehat{\Pi}_3(\lambda,\nu,k) &
\begin{array}{c}\text{relative}\\[-3pt]\text{error}\end{array}&
\begin{array}{c}\text{relative}\\[-3pt]\text{error} \\[-3pt]\text{bound}\end{array}\\[5pt]
\hline\hline\hline\hline\hline
         .9 & .5 & .2 & .50760 & .51315 & -.109\!\!\times\!\!10^{-1}& .867\!\!\times\!\!10^{-1}&.50770 &  -.212\!\!\times\!\!10^{-3}& .135\!\!\times\!\!10^{-2} \\
\hline .99 & .6 & .025 & .56514 & .56530 & -.287\!\!\times\!\!10^{-3}& .238\!\!\times\!\!10^{-2}& .56514& -.606\!\!\times\!\!10^{-7}& .403\!\!\times\!\!10^{-6}\\
\hline .999 & .75 & .004 & .60555 & .60556& -.682\!\!\times\!\!10^{-5} &.376\!\!\times\!\!10^{-4}& .60555 & -.138\!\!\times\!\!10^{-10}& .106\!\!\times\!\!10^{-9}\\
\hline .9999 & .8 & .0005 & .62453 & .62453& -.153\!\!\times\!\!10^{-6} &.110\!\!\times\!\!10^{-5}&.62453 & -.888\!\!\times\!\!10^{-15}& .597\!\!\times\!\!10^{-13} \\
\hline .999999 & .95 & .00002 & .71429 & .71429& .172\!\!\times\!\!10^{-8} & .540\!\!\times\!\!10^{-8}&.71429 & .190\!\!\times\!\!10^{-18}& .793\!\!\times\!\!10^{-18}\\
\hline\hline\hline\hline\hline
\end{array}
$$
} \textbf{Table 4.2.} \textsl{Numerical examples for the
approximation \emph{(\ref{eq:Pi-Gauss-reverse})} with
$\lambda\to{1}$, $k\to{1}$, $(1-\lambda)/(1-k)\to{0}$ and $\nu=7$.
The sixth and the ninth columns represents the relative errors
$R_{2,N}(\lambda,\nu,k)/\Pi(\lambda,\nu,k)$ in
\emph{(\ref{eq:Pi-Gauss-reverse})}. The seventh and the tenths
columns represent relative error bounds as given by the rhs of
\emph{(\ref{eq:R2-bound})-(\ref{eq:f})} divided by
$\Pi(\lambda,\nu,k)$.}

\paragraph{5. Acknowledgments.} We thank Professor Jos\'{e} Luis L\'{o}pez
for useful discussions and for reading the draft of this article.
The research of the first author has been supported by the Far
Eastern Branch of the Russian Academy of Sciences (grant
05-111-$\Gamma$-01-046) and Russian Basic Research Fund (grant
05-01-00099).


\begin{thebibliography}{99}
\bibitem{Carlson1} B.C.\,Carlson, Some series and bounds for incomplete elliptic integrals, \emph{J.
 Math. Phys.} {\bf 40}(1961), 125-134.
\bibitem{CarlsonBook} B.C.\,Carlson, \emph{Special Functions of Applied Mathematics}, Academic Press, New York, 1977.
\bibitem{Carlson2} B.C.\,Carlson and J.L.\,Gustafson,
Asymptotic approximations for symmetric elliptic integrals,
\emph{SIAM J. Math. Anal.} \textbf{25}(1994), 288-303.
\bibitem{Bat1} A.\,Erd\'{e}lyi, W.\,Magnus, F.\,Oberhettinger and
F.G.\,Tricomi, \emph{Higher transcendental functions, Vol. 1},
McGraw-Hill Book Company, Inc., New York, 1953.
\bibitem{Bat3} A.\,Erd\'{e}lyi, W.\,Magnus, F.\,Oberhettinger and
F.G.\,Tricomi,  \emph{Higher transcendental functions, Vol. 3},
McGraw-Hill Book Company, Inc., New York, 1955.
\bibitem{Lopez} J.L.\,L\'{o}pez, Asymptotic expansions of symmetric standard
elliptic integrals, \emph{SIAM J. Math. Anal.} \textbf{31},
4(2000), 754-775.
\bibitem{Lopez1} J.L.\,L\'{o}pez, Uniform asymptotic expansions of symmetric elliptic
integrals, \emph{Constructive Approximation} \textbf{17}, 4(2001),
535-559.
\bibitem{Radon} B.\,Radon, Sviluppi in serie degli integrali
ellipttici, \emph{Atti. Accad. Naz. Lincei, Mem., Cl. Sci. Fis.
Mat. Nat.} Ser.(8) \textbf{2}(1950), 69-109.
\bibitem{Riordan} J.\,Riordan, \emph{Combinatorial
Identities}, John Wiley and Sons, New York, 1968.
\bibitem{Prud3} A.P.\,Prudnikov,  Yu.A.\,Brychkov and
O.I.\,Marichev, \emph{Integrals and series, Volume 3: More Special
Functions}, Gordon and Breach Science Publishers, 1990.
\end{thebibliography}
\end{document}